\documentclass[12pt]{amsart}
\usepackage{amssymb}
\usepackage{amscd}

\newcommand{\N}{\mathbb{N}}
\newcommand{\C}{\mathbb{C}}

\newcommand{\Aa}{\mathbb{A}}

\newcommand{\Pp}{\mathbb{P}}

\newcommand{\Fq}{\mathbb{F}_q}

\newcommand{\Fqbar}{\overline{\mathbb{F}_q}}

\newcommand{\Ind}{\mathrm{Ind}}

\newcommand{\tr}{\mathrm{tr}}

\newcommand{\ch}{\mathrm{ch}}

\newcommand{\isomto}{\overset{\sim}{\rightarrow}}

\numberwithin{equation}{section}
\newtheorem{theorem}{Theorem}[section]
\newtheorem{corollary}[theorem]{Corollary}
\newtheorem{lemma}[theorem]{Lemma}
\newtheorem{proposition}[theorem]{Proposition}

\theoremstyle{definition}
\newtheorem{definition}[theorem]{Definition}

\newtheorem{example}[theorem]{Example}
\newtheorem{remark}[theorem]{Remark}
\newtheorem{problem}[theorem]{Problem}

\newcommand{\bdf}{\begin{definition}}
\newcommand{\edf}{\end{definition}\noindent}
\newcommand{\bex}{\begin{example}}
\newcommand{\eex}{\end{example}\noindent}
\newcommand{\bpr}{\begin{proposition}}
\newcommand{\epr}{\end{proposition}}
\newcommand{\blm}{\begin{lemma}}
\newcommand{\elm}{\end{lemma}}
\newcommand{\bth}{\begin{theorem}}
\renewcommand{\eth}{\end{theorem}}
\newcommand{\bpf}{\begin{proof}}
\newcommand{\epf}{\end{proof}\noindent}
\newcommand{\bcr}{\begin{corollary}}
\newcommand{\ecr}{\end{corollary}\noindent}
\newcommand{\beq}{\begin{equation}}
\newcommand{\eeq}{\end{equation}}
\newcommand{\bes}{\begin{equation*}}
\newcommand{\ees}{\end{equation*}}
\newcommand{\ben}{\begin{enumerate}}
\newcommand{\een}{\end{enumerate}}

\begin{document}
\title[Cohomology of the regular elements]
{The symmetric group representation on cohomology of the
regular elements of a maximal torus of the special linear group}
\author{Anthony Henderson}
\address{School of Mathematics and Statistics,
University of Sydney, NSW 2006, AUSTRALIA}
\email{anthonyh@maths.usyd.edu.au}
\thanks{This work was supported by Australian Research Council grant DP0344185}
%%%%%%%%%%%%%%%%%%%%%%%%%%%%%%%%%%%%%%%%%%%%%%%%%%%%%%%%
\begin{abstract}
We give a formula for the character of the representation of the
symmetric group $S_n$ on each isotypic component of the cohomology of the
set of regular elements
of a maximal torus of $SL_n$, with respect to the action of the centre.
\end{abstract}
%%%%%%%%%%%%%%%%%%%%%%%%%%%%%%%%%%%%%%%%%%%%%%%%%%%%%%%%%
\maketitle
%%%%%%%%%%%%%%%%%%%%%%%%%%%%%%%%%%%%%%%%%%%%%%%%%%%%%%%%%%%%%
\section{Introduction}
Let $n$ be a positive integer. Define the hyperplane complement
\begin{equation*}
T(1,n):=\{(z_1,z_2,\cdots,z_n)\in\C^n\,|\,
z_i\neq 0,\,\forall i,\, z_i\neq z_j,\,\forall i\neq j\}.
\end{equation*}
The symmetric group $S_n$ acts on $T(1,n)$ by permuting coordinates;
we can identify $T(1,n)$ with the set of regular elements of a maximal
torus of $GL_n(\C)$, and $S_n$ with the Weyl group of
the maximal torus.
This action induces representations of $S_n$ on the
cohomology groups $H^i(T(1,n))$ (taken with complex
coefficients). The characters of these representations are well known
(they will be encapsulated in a single `equivariant generating function'
in Theorem \ref{lehrerthm} below). For the purposes of this Introduction,
we recall only the `non-equivariant' information, i.e.\ the Betti
numbers of $T(1,n)$:
\beq 
\sum_i (-1)^i \dim H^i(T(1,n))\, q^{n-i}
= (q-1)(q-2)\cdots(q-n).
\eeq
These Betti numbers are particularly familiar, since $T(1,n)$ is
homotopy equivalent to the configuration space $C_{n+1}$ of 
$(n+1)$-tuples of distinct complex numbers.

Now consider the toral complement
\[ ST(1,n):=\{(z_1,z_2,\cdots,z_n)\in T(1,n)\,|\,
z_1z_2\cdots z_n=1\}. \]
This can be identified with the set of regular elements of a 
maximal torus
of $SL_n(\C)$. Of course $S_n$ still acts, and there is a commuting action
of $\mu_n$ (the centre of $SL_n(\C)$) by scaling. Thus we have a direct
sum decomposition of $S_n$-representations:
\begin{equation}
H^i(ST(1,n))\cong\bigoplus_{\chi\in\widehat{\mu_n}}
 H^i(ST(1,n))_\chi,
\end{equation}
where $H^i(ST(1,n))_\chi$ is the $\chi$-isotypic component of $H^i(ST(1,n))$.
In this paper we address the following problem, suggested by Lehrer:
\begin{problem} \label{mainprob}
Give a formula for the character of
the representation of $S_n$ on each $H^i(ST(1,n))_\chi$.
\end{problem}
\noindent
Our solution is given in \S3 (see especially \eqref{answereqn}). For now,
we state merely the non-equivariant version:
\beq \label{noneqeqn}
\begin{split}
&\sum_i (-1)^i \dim H^i(ST(1,n))_\chi\, q^{n-1-i}\\
&= \frac{(-1)^{n-n/r}n!}{r^{n/r}(\frac{n}{r})!}
(q-r-1)(q-2r-1)\cdots(q-(\frac{n}{r}-1)r-1),
\end{split}
\eeq
where $r$ is the order of $\chi$ in the character group $\widehat{\mu_n}$.

Summing over $\chi$, we deduce a formula \eqref{totaleqn}
for the character
of the total representation on $H^i(ST(1,n))$, of which the non-equivariant
specialization is:
\beq \label{noneqtotaleqn}
\begin{split}
&\sum_i (-1)^i \dim H^i(ST(1,n))\, q^{n-1-i}\\
&= \sum_{r\mid n}\frac{(-1)^{n-n/r}\phi(r)n!}{r^{n/r}(\frac{n}{r})!}
(q-r-1)(q-2r-1)\cdots(q-(\frac{n}{r}-1)r-1).
\end{split}
\eeq
This total formula was essentially already known. Since
$ST(1,n)$ is `minimally pure' in the sense of \cite{dimcalehrer},
\eqref{noneqtotaleqn} follows from the fact that the right-hand side
counts the number of $\Fq$-points of the variety $ST(1,n)$ for all
prime powers $q$ which are congruent to $1$ mod $n$. More generally,
\eqref{totaleqn} follows from a count of fixed points of twisted
Frobenius maps, which was the result \cite[Theorem 5.8]{fleischmannjaniszczak}
of Fleischmann and Janiszczak. See Remark \ref{fjrem} below.

Clearly the quotient of $ST(1,n)$ by $\mu_n$ can be identified with 
$\Pp T(1,n)$, the image of $T(1,n)$ in $\Pp^{n-1}(\C)$. So in
the case where $\chi$ is the trivial character, we are dealing with
\begin{equation} 
H^i(ST(1,n))^{\mu_n}\cong H^i(\Pp T(1,n)).
\end{equation}
In general, as \eqref{noneqeqn} suggests,
the representation $H^i(ST(1,n))_\chi$ is induced from the
wreath product subgroup $W(r,n/r):=\mu_r\wr S_{n/r}$
of $S_n$. To be precise, define the hyperplane complement
\[ T(r,m):= \{(z_1,z_2,\cdots,z_m)\in \C^{m}\,|\,
z_i\neq 0,\,\forall i,\, z_i^r\neq z_j^r,\,\forall i\neq j\}, \]
and its image $\Pp T(r,m)$ in $\Pp^{m-1}(\C)$. For $r\geq 2$,
we identify the corresponding
reflection group $G(r,1,m)$ with the wreath product $W(r,m)$
(if $r=1$, $W(1,m)=S_m$ acts on $T(1,m)$ as seen above).
In the following Theorem, $\varepsilon_n$ denotes the
sign character of $S_n$, and $\det_{n/r}$
the determinant character of $GL_{n/r}(\C)$, restricted to
$W(r,n/r)$.
\begin{theorem} \label{inductionthm}
Let $r$ be the order of $\chi\in\widehat{\mu_n}$. For every $i$, we 
have an isomorphism of representations of $S_n$:
\[ H^i(ST(1,n))_\chi\cong\varepsilon_n\otimes
\Ind_{W(r,n/r)}^{S_n}(\det{}_{n/r}\otimes 
H^{i-n+n/r}(\Pp T(r,n/r))). \]
\end{theorem}
\noindent
The proof of this Theorem given in \S4 merely equates the characters
of both sides; a more conceptual understanding of the isomorphism (or rather
the related isomorphism \eqref{mysteryeqn}), involving Orlik-Solomon-style 
bases for the cohomology groups, is given in \cite{mybases}.

When $\chi$ is faithful (i.e.\ $r=n$), Theorem \ref{inductionthm} 
says that
\begin{equation} \label{faithfuleqn}
H^i(ST(1,n))_\chi\cong
\left\{\begin{array}{cl}
\varepsilon_n\otimes\Ind_{\mu_n}^{S_n}(\psi),&\text{ if $i=n-1$,}\\
0,&\text{ otherwise,}
\end{array}\right.
\end{equation}
where $\mu_n$ is embedded in $S_n$ as the subgroup generated by
an $n$-cycle, and $\psi\in\widehat{\mu_n}$ is any faithful character 
(it doesn't matter which -- note also that tensoring with $\varepsilon_n$
makes no difference unless $n\equiv 2$ (mod $4$)).
This result for prime $n$ was proved in \cite[\S4.4]{dimcalehrer}.
%%%%%%%%%%%%%%%%%%%%%%%%%%%%%%%%%%%%%%%%%%%%%%%%%%%%%%%%%%%%%%
\section{Equivariant weight polynomials}
Suppose $X$ is an irreducible complex variety which
is \emph{minimally pure} in the sense that $H_c^i(X)$ is a pure Hodge structure
of weight $2i-2\dim X$ for all $i$ (see \cite{dimcalehrer}). 
Let $\Gamma$ be a finite group acting on $X$.
We define the \emph{equivariant weight polynomials} of this action by
\[ P(\gamma, X, q):=\sum_i (-1)^i\, \tr(\gamma,H_c^i(X))\, q^{i-\dim X}, \]
for all $\gamma\in\Gamma$, where $q$ is an indeterminate ($=t^2$ in the 
notation of \cite{dimcalehrer}).
We also define
\[ P^\Gamma(X,q):=\sum_i (-1)^i\, [H_c^i(X)]\, 
q^{i-\dim X} \in R(\Gamma)[q], \]
where $R(\Gamma)$ is the complexified representation ring of $\Gamma$.
If $\Delta$ is an abelian finite group acting on $X$ whose action commutes
with that of $\Gamma$, and $\chi$ is a character of $\Delta$, we define
\begin{equation*}
\begin{split} 
P(\gamma,\chi,X,q)&:=\sum_i (-1)^i\, \tr(\gamma,H_c^i(X)_\chi)\, q^{i-\dim X}\\
&=\sum_i (-1)^i \frac{1}{|\Delta|}\sum_{\delta\in\Delta}
\chi(\delta)^{-1} \tr((\gamma,\delta),H_c^i(X)) \, q^{i-\dim X}\\
&=\frac{1}{|\Delta|}\sum_{\delta\in\Delta}\chi(\delta)^{-1}
P((\gamma,\delta),X,q),
\end{split}
\end{equation*}
and similarly
\[ P^{\Gamma}(\chi,X,q):=\sum_i (-1)^i\, [H_c^i(X)_\chi]\, q^{i-\dim X}\in 
R(\Gamma)[q]. \]
If $X$ is nonsingular, we can translate knowledge of $H_c^i(X)$
and $H_c^i(X)_\chi$ into knowledge of $H^i(X)$ and $H^i(X)_\chi$ by
Poincar\'e duality. 

Now for any positive integers $r$ and $n$,
$T(r,n)$ (respectively, $\Pp T(r,n)$) 
is a nonsingular irreducible minimally pure variety of dimension $n$ 
(respectively, $n-1$); minimal purity is a 
standard property of hyperplane complements 
(\cite[Example 3.3]{dimcalehrer}).
Also, $ST(1,n)$ is clearly a nonsingular
irreducible variety of dimension $n-1$. To show that it is minimally pure,
one can use \cite[Corollary 4.2]{dimcalehrer}, or else observe that
it is the quotient of $\Pp T(n,n)$ by the free action of a finite group,
as follows.

Recall that $W(r,m)=S_m\rtimes\mu_r^m$ acts on $T(r,m)$ and $\Pp T(r,m)$;
the $S_m$ factor acts by permuting the coordinates, and $\mu_r^m$ acts
by scaling them.
Define a surjective map $\varphi:\Pp T(n,n)\to ST(1,n)$ by
\[ \varphi([x_1:x_2:\cdots:x_n])=\frac{1}{x_1 x_2\cdots x_n}(x_1^n,x_2^n,
\cdots,x_n^n). \]
The fibres of $\varphi$ are clearly the orbits of the normal
subgroup $S\mu_n^n\subset W(n,n)$ defined by
\[ S\mu_n^n:=\{(\zeta_1,\cdots,\zeta_n)\in\mu_n^n\,|\,
\zeta_1\zeta_2\cdots\zeta_n=1\}. \]
The action of $S\mu_n^n$ on $\Pp T(n,n)$ becomes free once one factors out the
subgroup $\{(\zeta,\zeta,\cdots,\zeta)\}$, which acts trivially.
Thus $ST(1,n)$ is minimally pure, and solving
Problem \ref{mainprob} amounts to computing the polynomials
$P(w,\chi,ST(1,n),q)$, for all $w\in S_n$ and 
$\chi\in\widehat{\mu_n}$.

Consider the quotient of $T(n,m)$ by $S\mu_n^m$ for arbitrary
$m\geq 1$. This can be identified with
\[ T^{(n)}(1,m)
:=\{((z_i),z)\in T(1,m)\times\C^\times \,|\, z^n=z_1\cdots z_m\}. \]
The quotient map $\psi:T(n,m)\to T^{(n)}(1,m)$ is given by
\[\psi(x_1,x_2,\cdots,x_m)=((x_1^n,x_2^n,\cdots,x_m^n),x_1 x_2\cdots x_m).\]
The group $S_m\times\mu_n\cong W(n,m)/S\mu_n^m$ acts on the quotient 
$T^{(n)}(1,m)$
in the obvious way: $S_m$ acts on the $T(1,m)$
component, and $\mu_n$ acts by scaling $z$. 

When $m=n$, we have an isomorphism 
\begin{equation*}
\begin{split}
T^{(n)}(1,n)&\isomto ST(1,n)\times\C^\times\\
((z_1,\cdots,z_n),z)&\mapsto ((z_1 z^{-1},\cdots,z_n z^{-1}),z),
\end{split}
\end{equation*}
which respects the $S_n$-action, and transforms the $\mu_n$-action on
$T^{(n)}(1,n)$ into the inverse of the $\mu_n$-action on $ST(1,n)$,
and a scaling action on $\C^\times$. Since the latter has no effect on 
cohomology,
\begin{equation} \label{bundleeqn} 
P(w,\chi,ST(1,n),q)
=\frac{1}{q-1}P(w,\chi^{-1},T^{(n)}(1,n),q).
\end{equation}
So we aim to compute $P(w,\chi^{-1},T^{(n)}(1,n),q)$; it turns out to be
convenient to compute the polynomials  $P(w,\chi^{-1},T^{(n)}(1,m),q)$
for all $m\geq 1$ and $w\in S_m$ together.

\begin{remark} \label{gusrem}
One can see \textit{a priori} that allowing
$m\neq n$ incurs no extra work, thanks to the following neat argument, 
pointed out to me by Lehrer. If $d=\gcd(m,n)$, the action of $\mu_{n/d}
\subset\mu_n$
on $T^{(n)}(1,m)$ is part of the action of the connected
group $\C^\times$ defined by
\[ t.((z_i),z)=((t^{n/d}z_i),t^{m/d}z). \]
Hence $\mu_{n/d}$ acts trivially on cohomology, 
so $P(w,\chi^{-1},T^{(n)}(1,m),q)=0$
unless $\chi|_{\mu_{n/d}}=1$, i.e.\ $\chi^d=1$, i.e.\ $r|m$,
where $r$ is the order of $\chi$. Moreover, if $r|m$, then
writing $\chi^\circ$ for the character of $\mu_r$ such that
$\chi(\zeta)=\chi^\circ(\zeta^{n/r})$ for all $\zeta\in\mu_n$,
and $\chi'$ for the character of $\mu_m$ defined by
$\chi'(\zeta)=\chi^\circ(\zeta^{m/r})$ for all $\zeta\in\mu_m$,
we have
\begin{equation*}
\begin{split}
P(w,\chi^{-1},T^{(n)}(1,m),q)&= P(w,(\chi^\circ)^{-1},T^{(r)}(1,m),q)\\
&= P(w,(\chi')^{-1},T^{(m)}(1,m),q).
\end{split}
\end{equation*}
We will not actually use this observation.
\end{remark}

The identification of $T^{(n)}(1,m)$ with the quotient of $T(n,m)$ by 
$S\mu_n^m$ has the following consequence for equivariant weight polynomials:
\begin{proposition} \label{sumprop}
For any $w\in S_m$ and $\chi\in\widehat{\mu_n}$,
\begin{equation*}
\begin{split} 
P(w,\chi^{-1},&T^{(n)}(1,m),q)\\
&=\frac{1}{n^m}
\sum_{(\zeta_i)\in\mu_n^m} \chi(\zeta_1\cdots\zeta_m)
P(w(\zeta_1,\cdots,\zeta_m),T(n,m),q).
\end{split}
\end{equation*}
\end{proposition}
\begin{proof}
It is well known that if $V$ is a representation of the finite group
$G$ and $V^{H}$ is the subspace invariant under the normal subgroup
$H$, the character of $G/H$ on $V^H$ is given by
\[ \tr(gH,V^H)=\frac{1}{|H|}\sum_{h\in H}\tr(gh, V). \]
Now apply this with $V=H_c^i(T(n,m))$, $G=W(n,m)$, and $H=S\mu_n^m$,
so that $V^H\cong H_c^i(T^{(n)}(1,m))$ and $G/H\cong S_m\times\mu_n$.
We find that for all $\zeta\in\mu_n$,
\begin{equation*}
P((w,\zeta),T^{(n)}(1,m),q)=\frac{1}{n^{m-1}}
\sum_{\substack{(\zeta_i)\in\mu_n^m\\
\zeta_1\cdots\zeta_m=\zeta}}
P(w(\zeta_1,\cdots,\zeta_m),T(n,m),q).
\end{equation*}
Combining this with the fact that
\begin{equation*}
P(w,\chi^{-1},T^{(n)}(1,m),q)=\frac{1}{n}\sum_{\zeta\in\mu_n} \chi(\zeta)
P((w,\zeta),T^{(n)}(1,m),q)
\end{equation*}
gives the desired result.
\end{proof}
%%%%%%%%%%%%%%%%%%%%%%%%%%%%%%%%%%%%%%%%%%%%%%%%%%%%%%%%%%%%%%
\section{Generating functions}
In this section we will
compute the sum in Proposition \ref{sumprop}
using the known formula for the equivariant
weight polynomials of $T(r,m)$. As is usual in dealing with
characters of symmetric groups and wreath products, the
computations become easier if one uses suitable `equivariant
generating functions'.

For any $r\geq 1$, let $\Lambda(r)$ denote the polynomial ring
$\C[p_i(\zeta)]$ in countably many independent variables $p_i(\zeta)$
where $i$ is a positive integer and $\zeta\in\mu_r$. Define an $\N$-grading
on $\Lambda(r)$ by $\deg(p_i(\zeta))=i$. 
Also let $\Lambda(r)[q]:=\Lambda(r)\otimes_\C
\C[q]$, with the $\N$-grading given by the first factor (so $\deg(q)=0$).
Let $\Aa(r)=\C[\![p_i(\zeta)]\!]$ be the completion of $\Lambda(r)$,
and set $\Aa(r)[q]=\Aa(r)\otimes\C[q]$.

As in \cite[Chapter I, Appendix B]{macdonald}, we define an isomorphism 
$\ch_{W(r,m)}:R(W(r,m))\isomto\Lambda(r)_m$ by
\[ \ch_{W(r,m)}([M])=\frac{1}{r^m m!}\sum_{y\in W(r,m)}\tr(y,M)\, p_y, \]
where $p_y=\prod_{i,\zeta}p_i(\zeta)^{a_i(\zeta)}$ if $y$ has
$a_i(\zeta)$ cycles of length $i$ and type $\zeta$. 
Note that to recover $\tr(y,M)$ from $\ch_{W(r,m)}([M])$ one must multiply the
coefficient of $\prod_{i,\zeta}p_i(\zeta)^{a_i(\zeta)}$ by the order of
the centralizer of $y$, which is 
$\prod_{i,\zeta}a_i(\zeta)!(ri)^{a_i(\zeta)}$.
Write $\ch_{W(r,m)}$ also for the induced
isomorphism $R(W(r,m))[q]\isomto\Lambda(r)[q]_m$.

The result we need on the equivariant weight polynomials
for $T(r,m)$ can be conveniently stated in terms of
the equivariant generating function $P(r,q)\in\Aa(r)[q]$ defined by
\begin{equation*}
\begin{split} 
P(r,q)&:=1+\sum_{m\geq 1}\ch_{W(r,m)}(P^{W(r,m)}(T(r,m),q))\\
&=1+\sum_{m\geq 1}\frac{1}{r^m m!}\sum_{y\in W(r,m)}P(y,T(r,m),q)\, p_y.
\end{split}
\end{equation*}
In the following result $\mu(d)$ denotes the M\"obius function.
\begin{theorem} \label{lehrerthm}
If $R_{r,i,\theta}:=\sum_{d|i}|\{\zeta\in\mu_r\,|\,\zeta^d=\theta\}|
\mu(d)(q^{i/d}-1)\,\in\C[q]$,
\[ P(r,q)=\prod_{\substack{i\geq 1\\\theta\in\mu_r}}(1+p_i(\theta))
^{R_{r,i,\theta}/ri}. \]
\end{theorem}
\begin{proof}
This follows from Hanlon's result \cite[Corollary 2.3]{hanlonwreath}
on the M\"obius functions of Dowling lattices. Within the reflection
group context, it follows from results of Lehrer
in \cite{lehrerone} ($r=1$),
\cite{hyperoctahedral} ($r=2$) and \cite{lehrertwo}, \cite{blairlehrer} 
(general $r$). A short proof for all $r\geq 2$, based on an `equivariant
inclusion-exclusion' argument of Getzler, is given in
\cite[Theorem 8.4]{mywreath}
($T(r,m)$ is the same as what is there called $M(r,m)$). The $r=1$
case can be proved by the same method (note that $T(1,m)$ is different from 
the variety $M(1,m)$ considered in \cite[Theorem 8.2]{mywreath},
since it has the extra condition of nonzero coordinates).
\end{proof}
\noindent
Recovering the traces of individual elements by the above rule, we get
an equivalent statement, closer to Hanlon's and Lehrer's:
if $y$ in $W(r,m)$
has $a_i(\zeta)$ cycles of length $i$ and type $\zeta$,
\begin{equation} \label{trmeqn}
P(y,T(r,m),q)=\prod_{\substack{i\geq 1\\\zeta\in\mu_r}}
R_{r,i,\zeta}(R_{r,i,\zeta}-ri)\cdots(R_{r,i,\zeta}-(a_i(\zeta)-1)ri).
\end{equation}
There is an alternative description of the polynomials $R_{r,i,\theta}$. Define
\[ R_i^{(r)}:=\sum_{\substack{d|i\\\gcd(d,r)=1}}\mu(d)(q^{i/d}-1)\,\in\C[q]. \]
\begin{lemma} \label{polylemma}
If $t(\theta)$ denotes the order of $\theta$,
\begin{equation*}
R_{r,i,\theta}=\negthickspace\negthickspace\negthickspace
\sum_{s | \gcd(r/t(\theta),i)}\negthickspace\negthickspace
s\mu(s)R_{i/s}^{(r)}.
\end{equation*}
\end{lemma}
\begin{proof}
Since $\mu_r$ is cyclic of order $r$, we have
\[ |\{\zeta\in\mu_r\,|\,\zeta^d=\theta\}|=\left\{
\begin{array}{cl}
\gcd(d,r),&\text{ if $\gcd(d,r)\,|\,(r/t(\theta))$,}\\
0,&\text{ otherwise.}
\end{array}\right. \]
Hence
\[ R_{r,i,\theta}=\negthickspace\negthickspace
\sum_{s | (r/t(\theta))}
s \negthickspace
\sum_{\substack{d|i\\\gcd(d,r)=s}} \negthickspace \mu(d) (q^{i/d}-1). \]
The sum over $d$ has no terms unless $s|i$, in which case it equals
\[ \sum_{\substack{d | (i/s)\\ \gcd(d,r)=1}}\negthickspace
\mu(ds)(q^{i/ds}-1). \]
Since $\gcd(d,r)=1$ implies $\mu(ds)=\mu(d)\mu(s)$, this is
$\mu(s)R_{i/s}^{(r)}$.
\end{proof}
\noindent
This Lemma makes it clear that the $r=2$ case
of \eqref{trmeqn} is indeed equivalent to 
\cite[Theorem 5.6]{hyperoctahedral}.

As for $\Pp T(r,m)$, we have that for all $y\in W(r,m)$,
\beq \label{ptrmeqn}
P(y,\Pp T(r,m),q)=\frac{1}{q-1}P(y,T(r,m),q).
\eeq
For this, one need only show that the isomorphism
\[ \varphi:T(r,m)\to\Pp T(r,m)\times\C^\times:
(z_1,\cdots,z_m)\mapsto ([z_1:\cdots:z_m],z_1) \]
induces a $W(r,m)$-equivariant map on cohomology.
It is enough to check that $w\circ\varphi$ and $\varphi\circ w$ are
homotopic for all $w$ in a set of generators for $W(r,m)$,
which is straightforward.

Now for any $\chi\in\widehat{\mu_n}$, define the generating function
\begin{equation*}
\begin{split}
P(\chi,q)&:=1+\sum_{m\geq 1} \ch_{S_m}(P^{S_m}(\chi^{-1},T^{(n)}(1,m),q))\\
&=1+\sum_{m\geq 1}\frac{1}{m!}\sum_{w\in S_m}
P(w,\chi^{-1},T^{(n)}(1,m),q)\, p_w\,\in\Aa(1)[q].
\end{split}
\end{equation*}
We want a formula for this similar to
Theorem \ref{lehrerthm}. Define
\begin{equation*}
P_i^{(r)}:=\negthickspace\prod_{s|\gcd(r,i)}\negthickspace
(1-(-p_i)^{r/s})^{s\mu(s)R_{i/s}^{(r)}/ri}\ \in\Aa(1)[q].
\end{equation*}
\begin{theorem} \label{formulathm}
If $\chi\in\widehat{\mu_n}$ has order $r$,
$P(\chi,q)=\prod_{i\geq 1}P_i^{(r)}$.
\end{theorem}
\begin{proof}
Using Proposition \ref{sumprop}, we see that $P(\chi,q)$ equals
\begin{equation*}
1+\sum_{m\geq 1}\frac{1}{n^m m!}\sum_{\substack{w\in S_m\\
(\zeta_i)\in\mu_n^m}} \chi(\zeta_1\cdots\zeta_m)
P(w(\zeta_1,\cdots,\zeta_m),T(n,m),q)\, p_w,
\end{equation*}
which is precisely the result of applying to $P(n,q)$ the specialization
$p_i(\theta)\to\chi(\theta)p_i$.
So by Theorem \ref{lehrerthm},
\begin{equation*}
\begin{split}
P(\chi,q)
&=\exp\sum_{\substack{i\geq 1\\\theta\in\mu_n}}
\frac{R_{n,i,\theta}}{ni}\log(1+\chi(\theta)p_i)\\
&=\exp\sum_{\substack{i\geq 1\\\theta\in\mu_n\\m\geq 1}}
\frac{-R_{n,i,\theta}}{nmi}\chi(\theta)^m (-p_i)^m\\
&=\exp\sum_{\substack{i\geq 1\\d|i\\\zeta\in\mu_n\\m\geq 1}}
\frac{-\mu(d)}{nmi}\chi(\zeta)^{dm}(q^{i/d}-1)(-p_i)^m\\
&=\exp\negthickspace\negthickspace\negthickspace
\sum_{\substack{i\geq 1\\d|i\\m\geq 1,\, r|dm}}
\negthickspace\negthickspace
\frac{-\mu(d)}{mi}(q^{i/d}-1)(-p_i)^m,
\end{split}
\end{equation*}
since $\sum_{\zeta\in\mu_n}\chi(\zeta)^{dm}=n$ if $\chi^{dm}$ is trivial,
and vanishes otherwise. Thus we need only show that
\begin{equation} \label{fjeqn}
P_i^{(r)}=\exp\negthickspace\negthickspace
\sum_{\substack{d|i\\m\geq 1,\, r|dm}}
\negthickspace\negthickspace
\frac{-\mu(d)}{mi}(q^{i/d}-1)(-p_i)^m.
\end{equation}
The condition $r|dm$ is equivalent to $(r/\gcd(d,r))|m$. Writing
$s$ for $\gcd(d,r)$, the right-hand side becomes
\begin{equation*}
\begin{split}
\exp&\negthickspace\negthickspace\negthickspace
\sum_{\substack{s|r\\d|i,\,\gcd(d,r)=s\\m\geq 1,\,
(r/s)|m}} \negthickspace\negthickspace
\frac{-\mu(d)}{mi}(q^{i/d}-1)(-p_{i})^m\\
&=\exp
\sum_{s | \gcd(r,i)} \frac{s}{ri}\,\log(1-(-p_{i})^{r/s})\negthickspace
\sum_{\substack{d | i\\\gcd(d,r)=s}}\negthickspace
\mu(d)(q^{i/d}-1).
\end{split}
\end{equation*}
By the same argument as in the proof of Lemma \ref{polylemma},
the sum over $d$ equals $\mu(s)R_{i/s}^{(r)}$ as required.
\end{proof}
\noindent
Note that $P(\chi,q)$ depends only on $r$, not on $n$ or $\chi$, and that,
as predicted in Remark \ref{gusrem}, 
its nonzero homogeneous components all have degree divisible by $r$.

There is no formula as neat as \eqref{trmeqn} for the individual
polynomials $P(w,\chi^{-1},T^{(n)}(1,m),q)$. But if $w\in S_m$ has $a_i$ cycles
of length $i$, we know that
\begin{equation} \label{chitn1meqn}
P(w,\chi^{-1},T^{(n)}(1,m),q)=\prod_{i\geq 1}a_i!\, i^{a_i}
(\text{coefficient of $p_i^{a_i}$ in $P_i^{(r)}$}).
\end{equation}
Note that for the right-hand side to be nonzero, $a_i$ must be divisible by
$(r/\gcd(r,i))$ for all $i$.
In the special case that $\gcd(r,i)=1$, we have
\[ P_i^{(r)}=(1-(-p_i)^r)^{R_i^{{(r)}}\!/ri}, \]
and the coefficient of $p_i^{a_i}$, where $a_i$ is divisible by $r$, is
\begin{equation}
(-1)^{a_i-a_i/r}\frac{R_i^{(r)}(R_i^{(r)}-ri)(R_i^{(r)}-2ri)\cdots
(R_i^{(r)}-(a_i-r)i)}{(ri)^{a_i/r}(a_i/r)!}.
\end{equation}
Some further special cases of note: the $r=1$ case of Theorem \ref{formulathm}
says that 
\[ P(\mathrm{triv},q)=\prod_{i\geq 1}(1+p_i)^{R_i^{(1)}\!/i}=P(1,q), \]
reflecting the fact that the quotient of $T^{(n)}(1,m)$ by $\mu_n$
is $T(1,m)$. Slightly more interesting is the $r=2$ case. We have
\[ P_i^{(2)}=\left\{\begin{array}{cl}
(1-p_i^2)^{R_i^{(2)}\!/2i},&\text{ if $i$ is odd,}\\
(1-p_i^2)^{R_i^{(2)}\!/2i}(1+p_i)^{-R_{i/2}^{(2)}/i},&\text{ if $i$ is even.}
\end{array}\right. \]
Hence if $i$ is even, the coefficient of $p_i^{a_i}$ is
\[ \sum_{j=0}^{\lfloor a_i/2\rfloor} (-1)^j\binom{R_i^{(2)}\!/2i}{j}
\binom{-R_{i/2}^{(2)}/i}{a_i-2j}. \]

Returning to Problem \ref{mainprob}, equations \eqref{bundleeqn} and
\eqref{chitn1meqn} tell us that if $w\in S_n$ has $a_i$ cycles of length
$i$,
\begin{equation} \label{answereqn}
P(w,\chi,ST(1,n),q)=\frac{1}{q-1}\prod_{i\geq 1}a_i!\, i^{a_i}
(\text{coefficient of $p_i^{a_i}$ in $P_i^{(r)}$}).
\end{equation}
Since there are $\phi(r)$ characters $\chi\in\widehat{\mu_n}$ of order $r$,
we deduce that
\begin{equation} \label{totaleqn}
P(w,ST(1,n),q)=\frac{1}{q-1}\sum_{r|n}\phi(r)\prod_{i\geq 1}a_i!\, i^{a_i}
(\text{coefficient of $p_i^{a_i}$ in $P_i^{(r)}$}).
\end{equation}
\begin{remark} \label{fjrem}
As mentioned in the Introduction,
if $q$ is specialized to a prime power congruent to 
$1$ mod $n$, the right-hand side of \eqref{totaleqn} equals the formula
given in
\cite[Theorem 5.8]{fleischmannjaniszczak} for the number of
$\mathbb{F}_q$-points of the regular set of a maximal torus of $SL_n(\Fq)$
obtained from a maximally split one by twisting with $w$.
(To see this, use the expression \eqref{fjeqn} for $P_i^{(r)}$; the
coefficient of $p_i^{a_i}$ is called $R_{a_i,n}^{i}(q)$ 
in \cite{fleischmannjaniszczak}.) 
We have effectively reproved that result, since
\[ P(w,ST(1,n),q)=|ST(1,n)(\Fqbar)^{wF}| \]
by Grothendieck's Frobenius trace formula and a suitable comparison theorem
of complex and $\ell$-adic cohomology; see 
\cite[5.3, Example 5.6]{dimcalehrer} for the details.
\end{remark}
%%%%%%%%%%%%%%%%%%%%%%%%%%%%%%%%%%%%%%%%%%%%%%%%%%%%%%%%%%%%%%%%
\section{Induction}
We now aim to prove Theorem \ref{inductionthm} by
interpreting the generating function $P(\chi,q)$
in terms of induced characters. Recall that $W(r,m)$ can be 
embedded in $S_{rm}$
as the centralizer of the product of $m$ disjoint $r$-cycles.
For any $\theta\in\mu_r$, let $t(\theta)$ denote the order of $\theta$.
\begin{lemma}
For any $W(r,m)$-module $M$,
\[ \ch_{S_{rm}}([\Ind_{W(r,m)}^{S_{rm}}(M)])=\ch_{W(r,m)}([M])|_{p_i(\theta)\to
p_{it(\theta)}^{r/t(\theta)}}. \]
\end{lemma}
\begin{proof}
This is a direct consequence of Frobenius' formula for induced characters,
once one observes that a cycle of length $i$ and type $\theta$ in
$W(r,m)$ beomes the product of $r/t(\theta)$ disjoint $it(\theta)$-cycles
when regarded as an element of $S_{rm}$.
\end{proof}
\begin{lemma} \label{indexlemma}
For any $W(r,m)$-module $M$,
\begin{equation*} 
\begin{split}
\ch_{S_{rm}}([\varepsilon_{rm}\otimes&\Ind_{W(r,m)}^{S_{rm}}
(\det{}^{-1}_{m}\otimes M)])\\
&=(-1)^{rm-m}\ch_{W(r,m)}([M])|_{p_i(\theta)\to
-\theta^{-1} (-p_{it(\theta)})^{r/t(\theta)}}.
\end{split}
\end{equation*}
\end{lemma}
\begin{proof}
If $y\in W(r,i)$ is a cycle of 
length $i$ and type $\theta$,
\[ \varepsilon_{ri}(y)\det{}_i(y)^{-1}
=(-1)^{i-1+(it(\theta)-1)r/t(\theta)}\theta^{-1}
=(-1)^{i-1+ri-r/t(\theta)}\theta^{-1}. \]
Also, if $y\in W(r,m)$,
\[ p_y|_{p_i(\theta)\to (-1)^{ri+i}p_i(\theta)}=(-1)^{rm-m}p_y. \]
Hence
\begin{equation*}
\begin{split} 
\ch_{W(r,m)}([\varepsilon_{rm}\otimes&\det{}^{-1}_{m}\otimes M])\\
&=(-1)^{rm-m}
\ch_{W(r,m)}([M])|_{p_i(\theta)\to
-\theta^{-1} (-1)^{r/t(\theta)}p_i(\theta)},
\end{split}
\end{equation*}
and the result follows by applying the previous Lemma.
\end{proof}

Now define an element $P'(r,q)\in\Aa(1)[q]$ by
\begin{equation*}
\begin{split}
&P'(r,q):=P(r,q)|_{p_i(\theta)\to
-\theta^{-1} (-p_{it(\theta)})^{r/t(\theta)}}\\
&=1+\sum_{m\geq 1}\ch_{W(r,m)}(P^{W(r,m)}(T(r,m),q))|_{p_i(\theta)\to
-\theta^{-1} (-p_{it(\theta)})^{r/t(\theta)}}\\
&=1+\sum_{m\geq 1}(-1)^{rm-m}\ch_{S_{rm}}(\varepsilon_{rm}\otimes
\Ind_{W(r,m)}^{S_{rm}}
(\det{}_{m}^{-1}\otimes P^{W(r,m)}(T(r,m),q))).
\end{split}
\end{equation*}
\begin{proposition} \label{equalityprop}
$P'(r,q)=\prod_{i\geq 1}P_i^{(r)}$.
\end{proposition}
\begin{proof}
By Theorem \ref{lehrerthm} and Lemma \ref{polylemma}, we have
\begin{equation*}
\begin{split}
P'(r,q)&=\prod_{\substack{i\geq 1\\\theta\in\mu_r}}
(1-\theta^{-1} (-p_{it(\theta)})^{r/t(\theta)})^{R_{r,i,\theta}/ri}\\
&=\negthickspace\negthickspace\negthickspace
\prod_{\substack{i\geq 1\\\theta\in\mu_r\\s|\gcd(r/t(\theta),i)}}\negthickspace
\negthickspace
(1-\theta^{-1} (-p_{it(\theta)})^{r/t(\theta)})^{s\mu(s)R_{i/s}^{(r)}/ri}.
\end{split}
\end{equation*}
Applying to this the M\"obius inversion formula for cyclotomic 
polynomials, in the form
\[ \prod_{\substack{\theta\in\mu_r\\t(\theta)=t}}(1-\theta^{-1} X)
=\prod_{u|t}(1-X^{t/u})^{\mu(u)}, \]
we obtain
\begin{equation*}
P'(r,q)=\negthickspace\negthickspace
\prod_{\substack{i\geq 1\\t|r\\s|\gcd(r/t,i)\\u|t}}\negthickspace\negthickspace
(1-(-p_{it})^{r/u})^{s\mu(s)\mu(u)R_{i/s}^{(r)}/ri}.
\end{equation*}
Write this as $\prod_{i\geq 1} Q_i^{(r)}$, where $Q_i^{(r)}$ is the product of
all factors involving the variable $p_i$. Thus
\begin{equation*}
\begin{split}
Q_i^{(r)}&=\exp\negthickspace\negthickspace\negthickspace\
\sum_{\substack{t|\gcd(r,i)\\s|\gcd(r/t,i/t)\\u|t}}
\negthickspace\negthickspace
\frac{st\mu(s)\mu(u)}{ri}R_{i/st}^{(r)}\log(1-(-p_i)^{r/u})\\
&=\exp\negthickspace
\sum_{\substack{v|\gcd(r,i)\\u|v\\s|(v/u)}}
\negthickspace
\frac{v\mu(s)\mu(u)}{ri}R_{i/v}^{(r)}\log(1-(-p_i)^{r/u}),
\end{split}
\end{equation*}
where we have set $v=st$. Since $\sum_{s|(v/u)}\mu(s)$ is nonzero if and only
if $u=v$, we find that $Q_i^{(r)}=P_i^{(r)}$ as required.
\end{proof}
\begin{corollary} \label{equalitycor}
If $\chi\in\widehat{\mu_n}$ has order $r$, $P(\chi,q)=P'(r,q)$.
\end{corollary}
\begin{proof}
Combine Theorem \ref{formulathm} and Proposition \ref{equalityprop}.
\end{proof}
\begin{corollary} \label{pcor}
If $\chi\in\widehat{\mu_n}$ has order $r$, and $r|m$, we have the following
equality in $R(S_m)[q]$:
\begin{equation*} 
\begin{split}
P^{S_{m}}&(\chi^{-1},T^{(n)}(1,m),q)\\
&=(-1)^{m-m/r}\varepsilon_{m}\otimes
\Ind_{W(r,m/r)}^{S_m}(\det{}^{-1}_{m/r}\otimes P^{W(r,m/r)}(T(r,m/r),q)).
\end{split}
\end{equation*}
\end{corollary}
\begin{proof}
Under the isomorphism $\ch_{S_m}$, the left-hand side corresponds to
the degree-$m$ term of $P(\chi,q)$, and the right-hand side corresponds to
the degree-$m$ term of $P'(r,q)$.
\end{proof}
\noindent
Taking coefficients of $q^{i-m}$ on both sides and multiplying by
$(-1)^i$, we get an isomorphism of $S_m$-modules:
\begin{equation} 
H_c^i(T^{(n)}(1,m))_{\chi^{-1}}\cong\varepsilon_{m}\otimes 
\Ind_{W(r,m/r)}^{S_m}(\det{}_{m/r}^{-1}\otimes H_c^{i-m+m/r}(T(r,m/r))).
\end{equation}
By Poincar\'e duality, this is equivalent to
\[ H^{2m-i}(T^{(n)}(1,m))_{\chi}\cong \varepsilon_{m}\otimes
\Ind_{W(r,m/r)}^{S_m}(\det{}_{m/r}\otimes H^{m/r+m-i}(T(r,m/r))), \]
which after replacing $2m-i$ by $i$ gives
\begin{equation} \label{mysteryeqn}
H^i(T^{(n)}(1,m))_{\chi}\cong \varepsilon_{m}\otimes
\Ind_{W(r,m/r)}^{S_m}(\det{}_{m/r}\otimes H^{i-m+m/r}(T(r,m/r))).
\end{equation}

Finally, we prove Theorem \ref{inductionthm}.
Equations \eqref{bundleeqn}, \eqref{ptrmeqn}, and
Corollary \ref{pcor} together imply
\begin{equation*}
\begin{split}
P^{S_n}&(\chi,ST(1,n),q)\\
&=\frac{(-1)^{n-n/r}}{q-1}\varepsilon_n\otimes
\Ind_{W(r,n/r)}^{S_n}(\det{}^{-1}_{n/r}\otimes P^{W(r,n/r)}(T(r,n/r),q))\\
&=(-1)^{n-n/r}\varepsilon_n\otimes
\Ind_{W(r,n/r)}^{S_n}(\det{}^{-1}_{n/r}\otimes P^{W(r,n/r)}(\Pp T(r,n/r),q)).
\end{split}
\end{equation*}
Taking coefficients of $q^{i-n+1}$ on both sides and multiplying by
$(-1)^i$, we get an isomorphism of $S_n$-modules:
\begin{equation} 
H_c^i(ST(1,n))_{\chi}\cong\varepsilon_n\otimes
\Ind_{W(r,n/r)}^{S_n}(\det{}^{-1}_{n/r}\otimes H_c^{i-n+n/r}(\Pp T(r,n/r))).
\end{equation}
Since the right-hand side depends only on $n$ and the order of $\chi$,
this remains true if $\chi$ is replaced by $\chi^{-1}$. Then
Theorem \ref{inductionthm} follows by Poincar\'e duality.
%%%%%%%%%%%%%%%%%%%%%%%%%%%%%%%%%%%%%%%%%%%%%%%%%%%%%%%%%%%%%%%
\bibliographystyle{siam}

\begin{thebibliography}{1}

%\bibitem{barcelo}
%{\sc H.~Barcelo}, {\em On the action of the symmetric group on the free {Lie}
%  algebra and the partition lattice}, J. Combin. Theory Ser. A, 55 (1990),
%  pp.~93--129.

\bibitem{blairlehrer}
{\sc J.~Blair and G.~I. Lehrer}, {\em Cohomology actions and centralisers
  in unitary reflection groups}, Proc. London Math. Soc. (3), 83 (2001),
  no.~3, pp.~582--604.

\bibitem{dimcalehrer}
{\sc A.~Dimca and G.~I. Lehrer}, {\em Purity and equivariant weight
  polynomials}, in Algebraic Groups and Lie Groups, vol.~9 of Austral. Math.
  Soc. Lect. Ser., Cambridge University Press, Cambridge, 1997, pp.~161--181.

%\bibitem{esv}
%{\sc H.~Esnault, V.~Schechtman, and E.~Viehweg}, {\em Cohomology of local
%  systems on the complement of hyperplanes}, Invent. Math., 109 (1992),
%  pp.~557--561.

\bibitem{fleischmannjaniszczak}
{\sc P.~Fleischmann and I.~Janiszczak}, {\em The number of regular semisimple
  elements for {Chevalley} groups of classical type}, J. Algebra, 155 (1993),
  pp.~482--528.

\bibitem{hanlonwreath}
{\sc P.~Hanlon}, {\em The characters of the wreath product group acting on
  the homology groups of the {Dowling} lattices}, J. Algebra, 91 (1984),
  pp.~430--463.

\bibitem{mywreath}
{\sc A.~Henderson}, {\em Representations of wreath products on cohomology of
  {De Concini-Procesi} compactifications}, Int. Math. Res. Not., 2004:20 
  (2004), pp.~983--1021.

\bibitem{mybases}
\leavevmode\vrule height 2pt depth -1.6pt width 23pt, {\em Bases for certain
  cohomology representations of the symmetric group}, to appear in
  J. Algebraic Combin., math.RT/0508162.

\bibitem{hyperoctahedral}
{\sc G.~I. Lehrer}, {\em On hyperoctahedral hyperplane complements}, in The
  Arcata Conference on Representations of Finite Groups (Arcata, Calif., 1986),
  vol.~47 of Proc. Sympos. Pure Math., Amer. Math. Soc., Providence, RI, 1987,
  pp.~219--234.

\bibitem{lehrerone}
\leavevmode\vrule height 2pt depth -1.6pt width 23pt, {\em On the {Poincar\'e}
  series associated with {Coxeter} group actions on complements of
  hyperplanes}, J. London Math. Soc. (2), 36 (1987), pp.~275--294.

\bibitem{lehrertwo}
\leavevmode\vrule height 2pt depth -1.6pt width 23pt, {\em {Poincar\'e}
  polynomials for unitary reflection groups}, Invent. Math., 120 (1995),
  pp.~411--425.

\bibitem{macdonald}
{\sc I.~G. Macdonald}, {\em Symmetric Functions and {Hall} Polynomials}, Oxford
  Univ. Press, second~ed., 1995.

\end{thebibliography}

%%%%%%%%%%%%%%%%%%%%%%%%%%%%%%%%%%%%%%%%%%%%%%%%%%%%%%%%%%%%%%%
\end{document}